\documentclass[english]{paper}
\usepackage{fullpage} 
\usepackage[T1]{fontenc}
\usepackage[latin9]{inputenc}
\usepackage{babel}

\usepackage{tipa}
\usepackage{amsthm}
\usepackage{amsmath,amsfonts,amssymb,amsbsy,amsthm}
\usepackage{graphicx}
\usepackage{wrapfig}
\usepackage{cite}
\makeatletter

\theoremstyle{plain}
\newtheorem{theorem}{Theorem}
\theoremstyle{definition}
\newtheorem{definition}{Definition}
\theoremstyle{plain}
\newtheorem{lemma}{Lemma}
\theoremstyle{plain}

\theoremstyle{plain}

\theoremstyle{plain}

\usepackage[nothing]{algorithm} \usepackage[noend]{algorithmic} 
\usepackage{float}
\newdimen\algorithmicindent \algorithmicindent=0.5cm

\usepackage{color}
\graphicspath{{figures/}{./}}

\newif\ifnotesw\noteswtrue

\newif\ifnotesw\noteswtrue

\def\blfootnote{\xdef\@thefnmark{}\@footnotetext}

\makeatother

\begin{document}

\title{Monotonic and Non-Monotonic Epidemiological Models on Networks}

\author{Alexander Gutfraind $^{1}$}

\date{\today}

\maketitle

\begin{abstract}
Contact networks can significantly change the course of epidemics,
affecting the rate of new infections and the mean size of an outbreak.
Despite this dependence, some characteristics of epidemics are not contingent on the contact network
and are probably predictable based only on the pathogen.
Here we consider SIR-like pathogens and give an elementary proof that 
for any network increasing the probability of transmission increases the mean outbreak size.
We also introduce a simple model, termed 2FleeSIR, in which susceptibles protect themselves by avoiding contacts with infectees.
The 2FleeSIR model is non-monotonic: for some networks, increasing transmissibility actually decreases the final extent.
The dynamics of 2FleeSIR are fundamentally different from SIR because 2FleeSIR exhibits no outbreak transition in densely-connected networks.
We show that in non-monotonic epidemics, public health officials might be able to intervene in a fundamentally new way
to change the network so as to control the effect of unexpectedly-high virulence. 
However, interventions that decrease transmissibility might actually cause more people to become infected.

\end{abstract}

\section{Introduction}
\setcounter{footnote}{2} \blfootnote{$^{1}$ Theoretical Division, Los Alamos National Laboratory, Los Alamos, New Mexico, 87545, USA, \href{mailto:agutfraind.research@gmail.com}{agutfraind.research@gmail.com}.}
Epidemic, infection, contagion and cascade are interchangeable
terms referring to a probabilistic process over a population connected
through a contact network. Such a process can model infectious diseases,
but is also relevant in studying phenomena such as transmission of information
and propagation of changes through a system~\cite{Newman03review,Volz08,Zager08,Coburn09}. 

In studying epidemics, one of the central problems is how the structure
of the network affects propagation.
For example, the epidemic might affect more people if infected individuals 
are connected to many susceptible individuals, rather than to other infectees~\cite{Meyers07}.
Another important problem is characterizing the final extent of the epidemic
as a function of the pathogen, such as its likelihood to transmit across contacts.

Intuition suggests that more transmissible epidemics and more dense contact
networks would lead to larger outbreaks on average.
We propose to call such epidemics
\emph{monotonic}.
In contrast, \emph{non-monotonic} epidemics are those that do not follow this rule.
With a non-monotonic epidemic, the extent might increase with transmissibility 
if the contact network has a particular structure but in other networks the extent might decrease.
We call the two possibilities concordant and discordant networks. 
Non-monotonic epidemics, if they exist, would be significant public-health challenge:
In them measures that retard transmission at the individual level (i.e. in a network of two connected individuals),
might create a concordant network and actually \emph{promote} the epidemic in the population at large.

This paper contributes an elementary rigorous proof showing that a certain 
commonly-used epidemic model is monotonic (sec.\ref{sec:monoton}).
More importantly it introduces a simple network-based model for behavior during an outbreak (sec.\ref{sec:non-mon}).
The latter model is shown as non-monotonic - a finding with both positive and negative implications to public health.

\section{Monotonicity of the SIR model\label{sec:monoton}}
Many epidemics, perhaps the majority, are described using the ``SIR'' model and its generalizations.
SIR models are useful in applications where the epidemic infects each network node at most once, for example
because it is lethal, or produces lasting immunity. 
Originally the SIR model was formulated as a system of three differential equations,
representing the number of {}``susceptible'' ($S)$, {}``infected''
($I$) and {}``removed''/''recovered'' ($R$) individuals~\cite{Ellner06}. 
An early network model of such epidemics was the percolation model
\cite{Kirkpatrick71} where individuals are sites of a D-dimensional discrete lattice.
A site in this lattice has probability $p$ of infecting a neighbor, independently for each neighbor \cite{Grassberger83}. 
Many generalizations of this model have been developed, for example, allowing for variability
in the contact rate between pairs, temporal changes in the network and others
\cite{Newman02,Volz08,Zager08,Segbroeck10,Lindquist11}. 

To make the discussion precise, let us consider a simple discrete-time SIR model on networks.
In this model, $\sigma$, $\tau$, and $\gamma$ are probability distributions
for the induction, transmission, and removal steps.
\begin{definition}
\label{def:Reed-Frost} The epidemic is a discrete-time
process on a directed graph $G(V,A)$. At each time $t=0,1,2\dots$
the set of nodes $V$ is partitioned into three disjoint sets: $V=S_{t}\cup I_{t}\cup R_{t}$.
At time $t=0$, $R_{0}=\emptyset$ and each node $u$ is added to the set of initial infectees, $I_{0}$,
with probability $\sigma(u)$ (the ``induction probability''). 
For each $t=0,1,2\dots$ the following
sequence of three events occurs for each node $u$ in state $I$ ($u\in I_{t}$):
(1) it is removed with probability $\gamma(u)$; (2) If $u$ was not
removed, then for each of its out-neighbors in state $S$, that is
node $v\in S_{t}$ where $(u,v)\in A$, there is independent probability
$\tau(u,v)$ of turning $v$ to state $I$ (giving $v\in I_{t+1}$);
(3) If $u$ is not removed yet, then it is {}``cured'' by time $t+1$
($u\in R_{t+1}$). Once in state $R$, the node remains there for
all future times. A node in state $S$ remains in that state unless
infected.
\end{definition}
The current definition is in some ways both more general and more narrow than those
considered in other works.  It is more general in that, for example, the set of initial
infectees is not necessarily a single node or even have a particular size
(although this is included as a special case by making $\sigma(v)\in\{0,1\}$ for all $v\in V$).
It is more narrow in that for instance, the time is discrete and the duration of every infection is exactly $1$.
A consequence of the discretization, which will be used shortly, is that 
the epidemic reaches a stationary state by time $t\leq n$,
where $n=\vert V \vert$. This is because any node infected
at time $t_{1}$ becomes $R$ at time $t_{1}+1$.  Thus at least one new
infection must be established at each $t$ or else the epidemic dies
out. 

How do the parameters of this epidemic affect its size?
Intuition suggests that increasing transmissibility $\tau$ would tend to increase the extent
of the epidemic. Heuristically when $\tau$ is larger
there more paths for the epidemic to reach any susceptible node. 
This argument is supported by analytic approximations,
results from special graphs and computer simulations on various variants of SIR \cite{Keeling99,Newman02,Volz08}.
However, increasing transmissibility sometimes introduces ``spread-blockers'' into the network: 
nodes that become infected early in the course of the epidemic but do not spread
it further, and then prevent the epidemic from spreading through them into other parts
of the network.  In certain network topologies, such an effect might conceivably
outweigh the increase in transmission probability.

We now show how the mean final extent of the epidemic, 
$\mathbb{E}\left(\left\vert R_{n} \right\vert \right)$, is affected by the parameters of the epidemic. 
Namely, the mean size of the epidemic on all contact network topologies
will not decrease because of (1) expansion in the set of initial infectees,
(2) an increase in the probability of transmission, or (3) decrease in the probability of removal:
\begin{theorem}
\label{thm:main}Let $\tau$, $\tau_{+}$ be two probability measures
for transmission events, with $\tau_{+}(e)\geq\tau(e)$ for all
$e\in A$. Let $\gamma_{-}$, $\gamma$ be two probability measures
for removal events, with $\gamma_{-}(v)\leq\gamma(v)$ for all
$v\in V$. Let $\sigma,\sigma_{+}\subset V$ be the induction probabilities into
the initial infectee sets of the two epidemics, with $\sigma_{+}(v)\geq \sigma(v)$ for all $v \in V$. 
Consider two epidemics
on $G$, epidemic $\overline{E}$ described by $\sigma_{+},\gamma_{-},\tau_{+}$
and epidemic $\underline{E}$ described by $\sigma,\gamma,\tau$. The mean
final extent $\mathbb{E}\left(\left\vert R_{n}\right\vert \right)$ is not smaller in $\overline{E}$:
\[
\mathbb{E}_{\sigma_{+},\gamma_{-},\tau_{+}}\left(\left\vert R_{n}\right\vert \right)\geq\mathbb{E}_{\sigma,\gamma,\tau}\left(\left\vert R_{n}\right\vert \right)\,.\]
\end{theorem}

The proofs in the appendix will show precisely when the equality occurs.  
A corollary of the monotonicity in $\tau$ is that adding edges (contacts)
to the network also never decreases the mean extent, because a missing edge $e$
is equivalent to $\tau(e)=0$.  Similarly, adding susceptible nodes cannot decrease
mean extent because of $\gamma$ monotonicity. 
Vaccinating nodes (in effect, increasing the probability of removal) would also tend
to decrease the final extent, never increasing it.

This theorem is an epidemiological equivalent of a much-studied
property of electrical networks, known as Rayleigh monotonicity~\cite{Doyle84,Choe08,Cibulka08}.
Rayleigh monotonicity describes the effect of individual resistors in a network
on the ``effective resistance'' of the entire network.
The effective resistance across any two terminals (nodes) imagines
that the network is just a single resistor and asks how large that resistor would have to be.
Rayleigh monotonicity is the finding that regardless of the network, 
increasing any single resistor in the network will not decrease the effective resistance
across any pair of terminals.
Interestingly, one of the many proofs of Rayleigh monotonicity is based on relating
it to a stochastic process on the network~\cite{Doyle84}.
Monotonicity results are also known for other stochastic contexts~\cite{Roerdink90,Haeggstroem99}.

Because SIR is such an important model, many previous studies considered the extent of SIR epidemics \cite{MartinLof86,picard90,Ball95,Ball10}. 
However, on networks most exact statements require $G$ to
have a special structure (e.g. a tree \cite{Newman02}). In spectral
graph theory some related results \cite{Draief06} show the connection
between the transmission probability $\tau$ ($\beta$ in their notation)
and cascade extent.
Most similar to this study is the superb work of Floyd, Kay and Shapiro~\cite{Floyd08}, 
who established the monotonicity of two SIR-type models using measure-theoretic arguments.
We offer an alternative proof here using elementary graph-theoretic methods and then introduce a novel non-monotonic model.
In outline, the proof uses graph cuts where a cut is formed from nodes that prevent the infection of a node $u$.

\section{Non-monotonic Models\label{sec:non-mon} and Self-Protective Behavior}
Can we expect monotonicity to hold in other epidemiological models?
Likely, variants of SIR as well as models with an exposed compartment (SEIR) should also exhibit monotonicity.
This is because increasing parameters such as transmission likelihood $\tau$ increases the effective
number of links in the network.

In light of this intuition, it is perhaps surprising that some models are non-monotonic.
Namely, there are epidemics where increasing the transmissibility actually decreases the extent.
The contact network may play a decisive role: while in some network topologies, parameters such as transmissibility
are positively related to the final extent, in other networks the relation is negative.
We will see that a simple source of non-monotonicity is the behavior of susceptible individuals who may attempt to reduce their risk of infection
- a process which was termed {\em reactive social distancing} \cite{Bootsma07}.
Depending on the pathogen they may obtain vaccinations, wear protective equipment, or maintain distance from fulminating individuals.

Such protective behaviors during epidemics have attracted a lot of recent research~\cite{Coburn09}, with models that describe 
the effect of infections on demand for vaccines, self-imposed isolation and other behaviors (see e.g. \cite{Epstein08,Funk09,Segbroeck10}.)
It is typical that individuals engage in protection only when they see the pathogen in a family member or a friend,
and usually dismiss news reports or advocacy by the authorities as exaggerated or too inconveniencing given the risk \cite{Glik07}.
One expects this lazy reaction particularly when the risk of infection is low or when the symptoms are rarely life-threatening, 
such as during seasonal respiratory infections.
Surprisingly, even during the severe 1918 Influenza pandemic people apparently delayed distancing until they witnessed fatalities \cite{Bootsma07}.
A lazy response is consistent with previous research that established that many costly behaviors are triggered
only after a social threshold is crossed, such as when more than one peer is already involved \cite{Centola07,Dodds04}.

Consider a simple model for this behavior that could be termed ``2FleeSIR''.  
This model combines the SIR model with a simplified model for distancing.
In 2FleeSIR the epidemic follows steps (1), (2) and (3) in Def.~\ref{def:Reed-Frost} of SIR.
2FleeSIR simply adds step (0) where any susceptible node $u$ that observes $\geq 2$ infected individuals 
in its immediate in-neighborhood will break all its ties and remain so for the duration of the epidemic.
Thus in effect, every edge $e$ coming into $u$ will have transmission $\tau(e)=0$. 
2FleeSIR is inspired by a model introduced in \cite{Funk09}
but supposes susceptibles base their behavior only on their immediate neighbors in the network.
2FleeSIR is a threshold cascade as in~\cite{Dodds04}, but unlike with other published models, here the contagion is harmful.
\begin{wrapfigure}{r}{0.4\columnwidth}%
\begin{centering}
\includegraphics[height=1in]{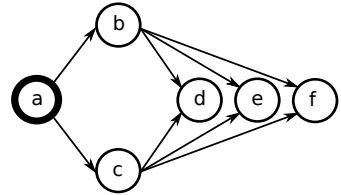}
\caption{\label{fig:nonmonotone}The ``Diamond'' graph illustrates the complex behavior of the 2FleeSIR model.  
For an epidemic that starts at $a$, increasing the transmission probability $\tau$ initially increases the epidemic's extent.
Then as $\tau\to 1$ the extent diminishes because individuals take protective measures.}
\par\end{centering}
\end{wrapfigure}
Thus when the threshold is crossed the propagation is suppressed rather than promoted.

\begin{figure}[!t]
\begin{centering}
\includegraphics[width=0.45\textwidth]{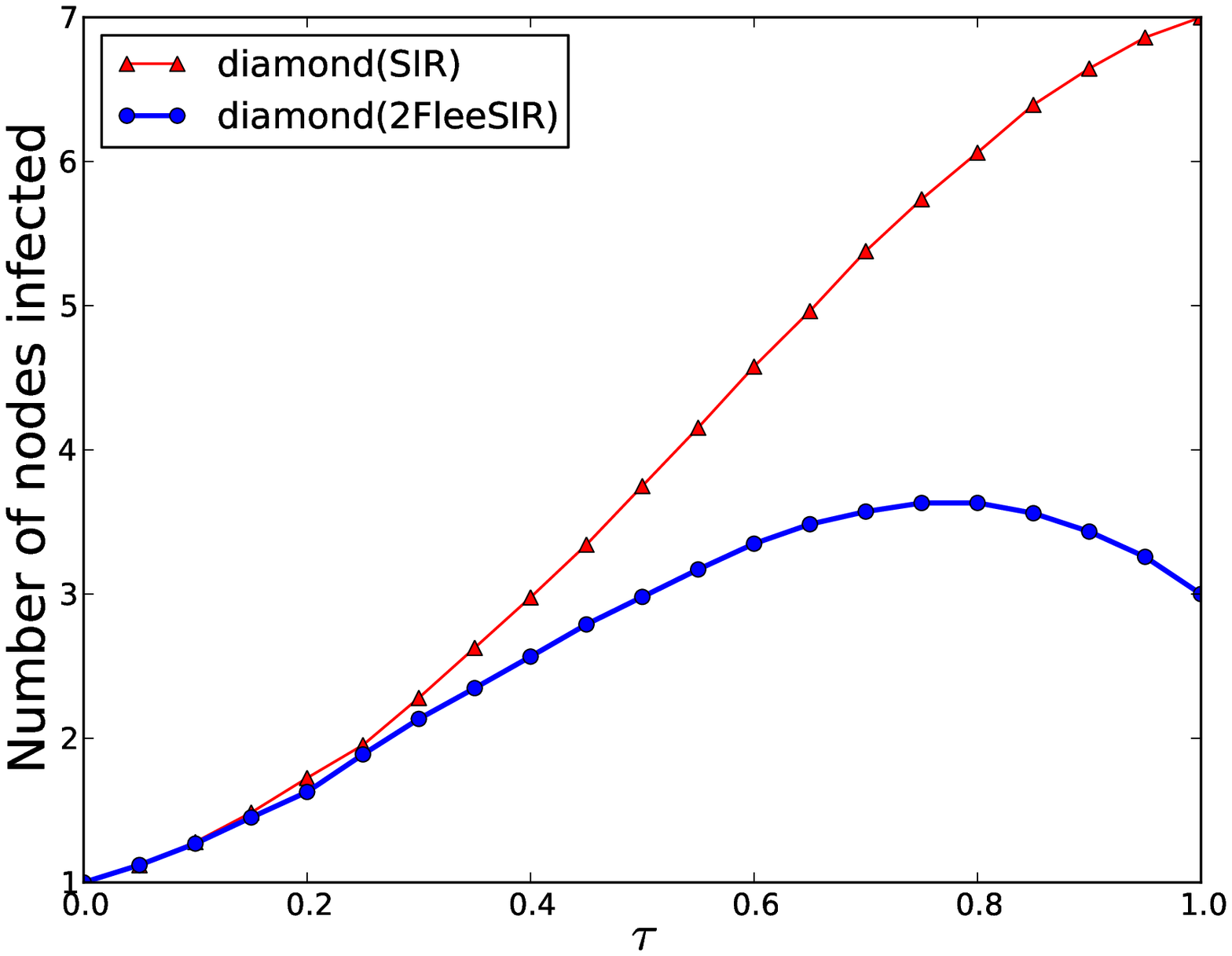}\hfill\includegraphics[width=0.45\textwidth]{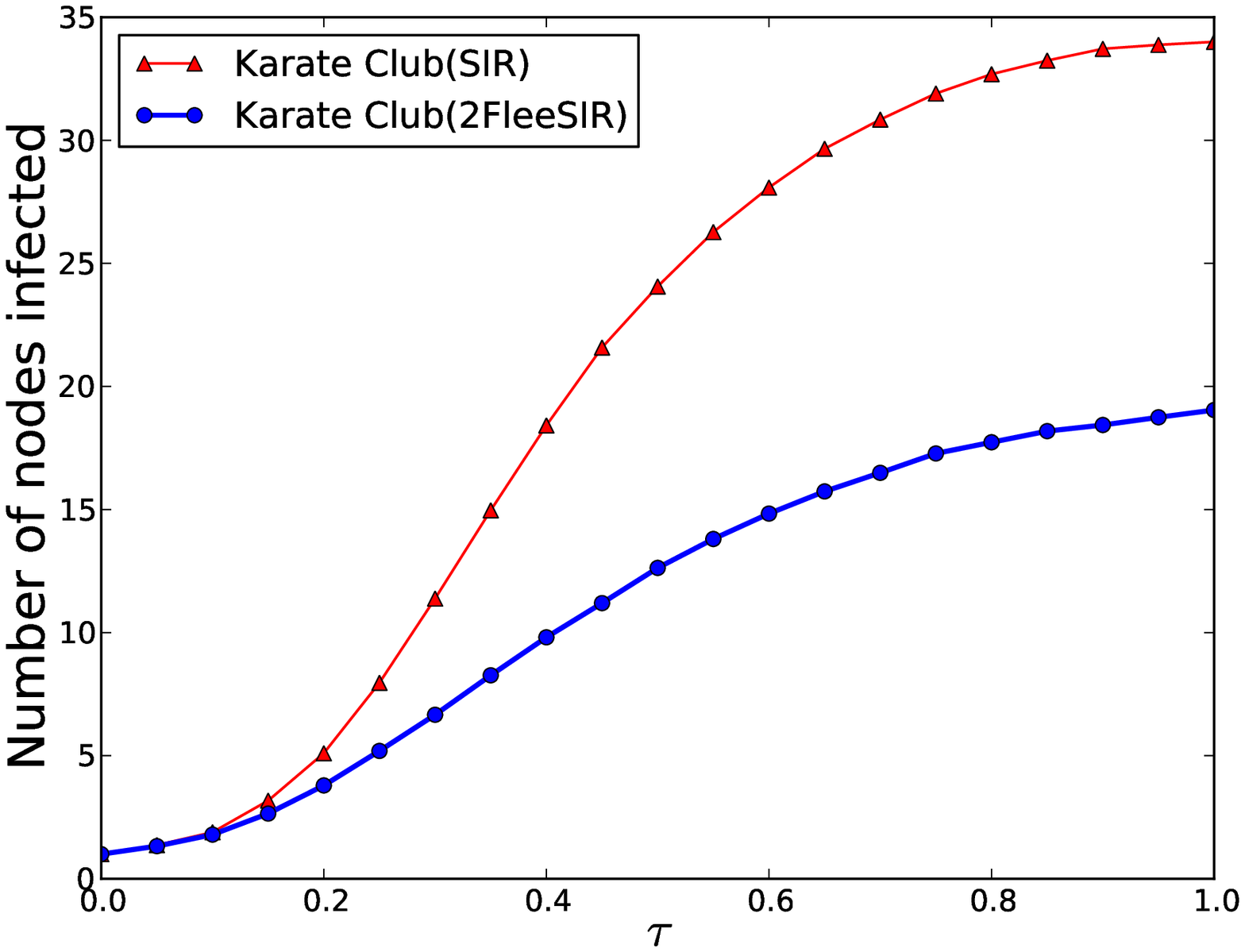}
\caption{\label{fig:nonmonotone-extent}The mean extent of SIR and 2FleeSIR epidemics on two contact networks:
on the Diamond network (left) and on an empirical social network, the Karate Club (right) \cite{Zachary77}.  
In both cases distancing in 2FleeSIR leads to a smaller mean extent as compared to SIR.  
In the Diamond network the extent actually decreases for sufficiently large $\tau$.
The experiments were slightly different: in the Diamond network the epidemic was started only at node $a$, 
while in the Karate Club it was started uniformly at all nodes. 
In the Diamond network, the decrease is only seen in epidemics that start at $a$.}
\par\end{centering}
\end{figure}

Consider now the spread of a 2FleeSIR epidemic on a small diamond-shaped network in Fig.~\ref{fig:nonmonotone}.
For simplicity, suppose that the epidemic originates at node $a$, the removal probability is $0$, and all edges transmit with probability $\tau$.
As $\tau$ increases it becomes increasingly likely that the epidemic infects both $b$ and $c$, which causes $d$,$e$ and $f$ to take protective measures.
As a result, its extent declines (see Fig.~\ref{fig:nonmonotone-extent}, left.)
It is difficult to determine if such a decline is seen frequently in empirical
contact networks because of insufficient data.
In simulations on other network we observed that it is more typical to see a reduction in the rate of increase in extent as a function of transmissibility
(see one simulation in Fig.~\ref{fig:nonmonotone-extent}, right).

There is typically much uncertainty as to parameters such as the basic reproductive ratio, $R_0$, of an epidemic~\cite{Meyers05}.  
Fortunately, if distancing is likely to happen, such a high $R_0$ is unlikely to produce dramatically many new cases
and may even reduce the total case load.
A striking effect of distancing is seen in the case of a fully-connected network. 
Recall that in fully-connected populations, the SIR model shows a jump in mean extent once $\tau$ passes the percolation transition $\tau=\frac{1}{n}$ \cite{Newman02}.
This outbreak transition is a key concern in attempting to control infections \cite{Ellner06,Meyers05}.
However, no such jump occurs in 2FleeSIR.
Rather when the epidemic starts at a single infected node, the ultimate extent is $1+\tau (n-1)$, that is, linear in $\tau$ and the number of nodes $n$.
This finding suggests that models like SIR might fundamentally misestimate the progression of epidemics if those cause wide-spread distancing behavior.
Sadly, in many epidemics the threshold for protective action might be difficult to cross 
or susceptibles might be unable to reliably identify infectives in their social neighborhoods.

The 2FleeSIR model also cautions about an interesting adverse effect of treatments.
Consider again the Diamond network, and suppose a vaccine or therapy completely prevents the infection of treated nodes
and suppose node $c$ is treated.
Under low $\tau$ it would reduce the final extent, but under $\tau\to 1$,
it would be impossible for node $d$, $e$ and $f$ to respond and they would all become infected.
As a result, the treatment would actually \emph{increase} the number of infectees from $3$ before treatment to $6$ after!
It is difficult to determine how frequently this is occurring in practice.

Earlier, Floyd et al. \cite{Floyd08}, pointed out that vaccination might interact with a non-monotonic epidemic in another way to cause an increase of the epidemic's extent.
Suppose the population is exposed to the epidemic in two phases such that in the second phase the pathogen has mutated into a more transmissible form.
Without vaccination there might be an early outbreak that creates recovered individuals
who act as cascade blockers.  With vaccination those blockers do not emerge
and the more virulent second strain is able to infect a larger number of individuals.

\section{Discussion}
Like this study, many recent projects in mathematical epidemiology developed network-based models of contagions,
often finding corrections to earlier models based on ordinary differential equations (ODEs).
But Theorem \ref{thm:main} here and Ref.\cite{Floyd08} show that in a fundamental sense the ODE version of SIR gives a correct description of the contagion: it is monotonic~\cite{Ellner06}.
In this narrow sense \emph{the network does not matter}.
We do not know how general this is, and an important avenue research is studying monotonicity in other epidemic models
from the SIR family, such as those with continuous time.  
Analogous results likely also hold in SIS, SEIR and other classes of models.

One may wonder whether monotonicity of final size might be true on average:
For example, given two epidemics, if $\overline{E}$ has a higher transmissibility on average (computed over edges),
would it also have a larger mean extent?  The answer is negative
because of the preponderance of network topology.
If the edges of lower transmission are bridge edges between parts of the contact network,
the effect of lower transmission would be magnified arbitrarily.

The 2FleeSIR model, introduced here, exemplifies the other possibility, rarely considered:
a non-monotonic epidemic where a more virulent epidemic is less likely to spread in certain discordant networks.
This finding has public health implications.
Discordant networks act like an automatic stabilizer, reducing the mean outbreak size if the virulence increases.
But, even if the pre-outbreak network is concordant, public health officials might be able to intervene 
and make the network discordant and thus passively control the outbreak.
Non-monotonicity also suggests that it is possible in theory for therapies such as vaccines
to increase the extent of epidemics by interfering with socially-driven preventative behavior of individuals.
At this point there is no hard data for speculation on whether such an effect occurs frequently in practice.
Indeed the lesson from recent outbreaks of Measles is just how dangerous insufficient vaccine coverage can be to public health.

\section*{Acknowledgements}
Milan Bradonji$\text{\textipa{\'{c}}}$ and Vadas Gintautas helped with critical suggestions and encouragement. Andrea Pugliese and Eli Ben-Naim suggested improvements. Part of this work was funded by the Department of Energy at the Los Alamos National Laboratory under contract DE-AC52-06NA25396
through the Laboratory Directed Research and Development program,
and by the Defense Threat Reduction Agency. Released as Los Alamos Unclassified Report 11-01481.

\appendix
\section{Elementary Proof of Monotonicity of SIR\label{sec:proofs}}
This applies graph-theoretical techniques to prove the monotonicity of SIR.
The proof builds on a Lemma: decreasing the removal rate increases the probability
of infecting every node in the network.
\begin{lemma}
\label{lem:gamma} Suppose two epidemics both originate at some source
node $s$, $\{s\}=I_{0}$ (i.e. $\sigma_{+}(v)=1=\sigma(v)$ iff $v=s$ otherwise $=0$), and both have $\tau(e)=1$ for all $e\in A$.
Let the probabilities of removal be $\gamma_{-}(v)$ and $\gamma(v)$
such $\gamma_{-}(v)\leq\gamma(v)$ for all $v\in V$. Then for any
node $u$, the probability of infection is not smaller under $\gamma_{-}$:\[
\mathbb{P}_{\gamma_{-}}\left(u\in R_{n}\right)\geq\mathbb{P}_{\gamma}\left(u\in R_{n}\right)\,.\]
\end{lemma}
\subsection{Proof of Lemma \ref{lem:gamma}}

\subsubsection*{Step 1: Special case}

Here we will prove the Lemma in the case where there is just one node
such that $\gamma_{-}(z)<\gamma(z)$, while $\gamma_{-}(v)=\gamma(v)$
for all $z\in V\smallsetminus z$.

Recall the distinction between removal events as opposed to {}``cure''
events: the former occur before infection of out-neighbors while the
latter after. Note now that it is possible to decide in advance of the
epidemic whether the event of removing $v$ occurs, conditional on
$v$ becoming infected (an example of the principle of deferred decision).
In this perspective, even if the epidemic has not yet happened, we
can speak of the probability that node $v$ is removed. 

The infection of $u$ at some time $k+1<n$ can occur if and only
if (1) there is at least one directed path from node $s$ to $u$: $(s,v_{1}),(v_{1},v_{2}),\dots,(v_{k},u)$
satisfying (2) $s\in I_{0},v_{1}\in I_{1},\dots,v_{k}\in I_{k}$,
and (3) where no $v_{i}$ is removed.  Because $\tau(e)=1$ for all $e\in A$,
if (3) is satisfied for path $p$ so does (2) for all nodes in $p$. 
Observe that any non-self-intersecting directed path (i.e. any path in a graph-theoretic sense)
from $s$ to $u$ meets the topological condition (1), and any infection
path must be a path $s\to u$ (since no node can be infected more than once). 

\begin{figure}[!ht]
\begin{centering}
\includegraphics[width=0.5\columnwidth]{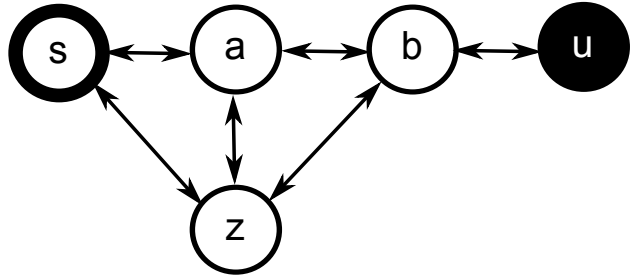}\caption{\label{fig:cuts}An illustration of $z$-vital cuts. For an epidemic that originates in node $s$, cut $\left\{ a,z\right\} $ is $s-u$ $z$-vital but cuts $\left\{ b\right\} $
and $\left\{ b,z\right\} $ are not.}

\par\end{centering}
\end{figure}

Let $PT$ be the set of all paths $s\to u$.
Assume for now that $PT\neq \emptyset$, i.e. condition (1) is satisfied for at least one path.
If for all $p\in PT$ at least one node $v\in p$ is removed (i.e. there
is an $s-u$ node cut in the graph) then $u$ cannot become infected.
Conversely, if $u$ is not infected then all of the paths have at
least one removed node.  Let $\mathbb{C}$ be the set of all possible
$s-u$ cuts in $(V,A)$.  
Recall that an $s-u$ cut $C$ is defined as a set of nodes so that any path $s\to u$ must pass a node $v\in C$.
Let $\mathbb{C}_{z}$ be the set of all {}``$z$-vital'' cuts: $\mathbb{C}_{z}$ is defined as all cuts
$C\in\mathbb{C}$ such that (1) $z\in C$ and, (2) if $z$ is excised 
from $C$, then $C$ is no longer an $s-u$ cut (illustrated in Fig.~\ref{fig:cuts}).

The event where $u$ is not infected could then be due to two types
of cuts: (a) $F_{z}=$the removed nodes form a cut $\in\mathbb{C}_{z}$
or (b) $F_{\bar{z}}=$the removed nodes form a cut $\in\mathbb{C}\smallsetminus\mathbb{C}_{z}$.
These two events are disjoint. We obtain:
\begin{eqnarray*}
\mathbb{P}(u\in R_{n}) & = & 1-\mathbb{P}(\mbox{a cut exists in the paths }s\to u)\\
 & = & 1-\left[\mathbb{P}(F_{z})+\mathbb{P}(F_{\bar{z}})\right].\end{eqnarray*}

Observe now that $\mathbb{P}_{\gamma}(F_{z})=\mathbb{P}_{\gamma}(F_{z}|z\mbox{ is removed})\underbrace{\mathbb{P_{\gamma}}(z\mbox{ is removed})}_{\gamma(z)}$.
The left term in this product does not depend on $\gamma(z)$ because
it expresses the probability of a type of $s-u$ cut 
in some graph $G'$ - a graph that in effect does not even have node $z$.
Observe also that the probability of the other event, $\mathbb{P}_{\tau}(F_{\bar{z}})$
does not depend on $\gamma(z)$ because any cut $\in F_{\bar{z}}$
will be a cut with or without the failure of $z$. Replacing terms
that do not depend on $\gamma$ with constants $c_{1}$ and $c_{2}$ get:
\begin{eqnarray*}
\mathbb{P}_{\gamma_{-}}\left(u\in R_{n}\right) & = & 1-\left[\mathbb{P}_{\gamma_{-}}(F_{z})+\mathbb{P}_{\gamma_{-}}(F_{\bar{z}})\right]\\
 & = & 1-c_{1}\gamma_{-}(z)-c_{2}\\
 & \geq & 1-c_{1}\gamma(z)-c_{2}  = \mathbb{P}_{\gamma}\left(u\in R_{n}\right)\,.
\end{eqnarray*}

This completes Step $1$. Observe that if $\gamma_{-}(z) < \gamma(z)$ (strictly), then
the probabilities of infection are equal in the two epidemics if and only if $c_{1}=0$, that is if $\mathbb{P}_{\gamma_{-}}(F_{z})=0$. 
This is the situation where no $z$-vital $s-u$ cut is possible or has positive probability. 
This means that either no directed path $s\to u$ exists 
(for example, $z$ might be a node of out-degree $0$)
or that for every $z$-vital cut $C$ the probability is $0$ that every node $t\neq z$ is removed.

\subsubsection*{Step 2: The general case}

Consider a sequence of epidemics $\overline{E}=E_{1},E_{2},\dots,E_{n}=\underline{E}$
where epidemic $i$ differs from epidemic $i+1$ at most at a single node
$v_{i}$. Namely, in epidemic $i$, $\gamma_{i}(u_{i})=\gamma_{-}(u_{i})$
while in epidemic $i+1$, $\gamma_{i+1}(u_{i})=\gamma(u_{i})\geq\gamma_{-}(u_{i})$.
Applying the result in Step 1 to each pair of epidemics produces a
sequence of inequalities $\mathbb{P}_{\gamma_{1}}\left(u\in R_{n}\right)\geq\mathbb{P}_{\gamma_{2}}\left(u\in R_{n}\right)\geq\dots\geq\mathbb{P}_{\gamma_{n}}\left(u\in R_{n}\right)$.
\qed

\subsection{Proof of Theorem \ref{thm:main} }

Recall that $\overline{E}$ is described by $\sigma_{+},\gamma_{-},\tau_{+}$
and $\underline{E}$ is described by $\sigma,\gamma,\tau$. The following
argument proves the result by constructing increasingly general situations.

\subsubsection*{Step 1: Simplest case}

Consider first the case where only the removal probabilities are different, namely: (1) $\overline{E}$ and $\underline{E}$
both originate in node $s$, and (2) in both $\tau(e)=1$ for
all edges $e$. By linearity of expectation and Lemma \ref{lem:gamma}:
\begin{eqnarray*}
\mathbb{E}_{\gamma_{-}}\left(\left\vert R_{n}\right\vert \right) & = & \sum_{u\in V}\mathbb{E}_{\gamma_{-}}\left({1}_{u\in R_{n}}\right) =  \sum_{u\in V}\mathbb{P}_{\gamma_{-}}\left(u\in R_{n}\right)\\
 & \geq & \sum_{u\in V}\mathbb{P}_{\gamma}\left(u\in R_{n}\right) = \mathbb{E}_{\gamma}\left(\left\vert R_{n}\right\vert \right)\,.
 \end{eqnarray*}

\subsubsection*{Step 2: Reduction of transmission events to removal events}

Any pair $\overline{E}$ and $\underline{E}$ which have possibly transmission probabilities different from $1$
can be mapped to a pair of the same extents but where the transmission probability is $1$, as follows.
Take epidemic $E_{\gamma,\tau}$ that starts at node $s$ and
has transmission probabilities $\tau$ and removal probabilities $\gamma$,
and construct an epidemic $E_{\hat{\gamma}}$ where $\tau(e)=1$ for
all edges $e$. This can be done by constructing graph $(V'',A'')$
in which every edge $e=(v,w)\in A$ is replaced by a two-hop path
$(v,h),(h,w)$ such that $\tau(v,h)=1=\tau(h,v)$ and the new helper
node $h$ has removal probability $\hat{\gamma}(h)=1-\tau(e)$. The
removal probability is unchanged for all regular nodes: $\hat{\gamma}(v)=\gamma(v)$
for all $v\in V$.

Couple the epidemics (any event $X$ in $E_{\gamma,\tau}$ occurs if and only if the corresponding event
occurs in $E_{\hat{\gamma}}$) and compare:
\begin{enumerate}
\item If in epidemic $E_{\gamma,\tau}$ node $v$ is infected then the infection
is transmitted to $w$ with probability $\tau(v,w)$. If in epidemic
$E_{\hat{\gamma}}$ node $v$ is infected then the infection is transmitted
to $w$ with the exact same probability: $$1\cdot\left[1-(1-\tau(v,w))\right]\cdot1=\tau(v,w)\,.$$
\item In epidemic $E_{\gamma,\tau}$ node $v$ is infected at time $t$
(transmission through $t$ edges) if and only if in epidemic $E_{\hat{\gamma}}$
node $v$ is infected at time $2t$ (because it takes two time steps
to transit every path between nodes copied from the original graph.
Consequently, in $E_{\hat{\gamma}}$ the epidemic will end by time
$2n$.
\item If $v\in V$ is removed in the final state of $E_{\gamma,\tau}$,
then the corresponding node is removed in $E_{\hat{\gamma}}$, and
vice versa.
\end{enumerate}
Therefore the final extents will be the same, correcting for removed
helper nodes (set $H$):\[
\mathbb{E}_{\hat{\gamma}}(\left\vert \widehat{R_{2n}}\smallsetminus H\right\vert )=\mathbb{E}_{\gamma,\tau}(\left\vert R_{n}\right\vert )\,.\]

Observe that when this reduction step is applied to $\overline{E}$
and $\underline{E}$, if for some edge $e$, $\tau_{+}(e)\geq\tau(e)$
then in the reduced epidemics $\hat{\gamma}_{-}\leq\hat{\gamma}$ in
the corresponding helper nodes. This satisfies the hypotheses of the theorem
and the restrictions of Step 1 ($\gamma_{-}(v)\leq\gamma(v)$ for all $v\in V$).

\subsubsection*{Step 3: Reduction of $\sigma$ and $\sigma_{+}$ to a single node $\alpha$}

Any pair $\overline{E}$ and $\underline{E}$ where there is more than one initial infectee
can be mapped to a pair of the same respective extents but where only one node is initially infected, as follows.

Take an epidemic $E_{\sigma,\gamma,\tau}$ that starts at time $t=0$ from an initial set formed with probabilities $\sigma$ on graph $(V,A)$ 
and spreads with transmission probabilities $\tau$ and removal probabilities $\gamma$. Couple it to an epidemic $E_{\alpha}$ that starts at time $t=-1$ at
a single node $\alpha$ in a slightly larger graph $(V',A')$:
\begin{enumerate}
\item Let $V'=V\cup\{\alpha\}$ and $I_{-1}=\{\alpha\}$
\item Let $\alpha$ never be removed before infecting neighbors: $\gamma(\alpha)=0$
\item For all $v\in V$ add an edge $(\alpha,v)$ : $A'=A\cup\{(\alpha,v)\mbox{ for all }v\in V\}$
\item The probability of transmission is $\tau(u,v)=\begin{cases}
\tau(u,v) & (u,v)\in A\\
\sigma(v) & u=\alpha \end{cases}\,.$
\end{enumerate}
The constructed epidemic originates at $\alpha$, and at time $t=0$ the set $I_{0}$ has
the same distribution as in the original epidemic.
The mean extent of $E_{\alpha}$ is just one node larger (node
$\alpha$) than the mean extent of $E_{\sigma,\gamma,\tau}$. This implies that
a situation where $\sigma_{+} \geq \sigma$ can be mapped to a problem
covered by Step 2, where both $\overline{E}$ and $\underline{E}$
originate at $\alpha$, because a larger $\sigma$ is equivalent to greater
probability of transmission $\tau_{+}$ for edges from $\alpha$ to
$q$ satisfying $\sigma_{+}(q) > \sigma(q)$.
This completes the theorem.\qed
An immediate consequence of Steps 2 and 3, and Lemma \ref{lem:gamma} is that every node is at least as likely to become infected
under $\overline{E}$ as under $\underline{E}$.

\bibliographystyle{IEEEtran}
\bibliography{netstorm,cascades}

\end{document}